\documentclass{amsart}
\usepackage{amssymb,amsmath}

\newcommand{\RR}{{\mathbb R}}

\newcommand{\e}{\varepsilon}

\newcommand{\Lap}{\Delta}
\newcommand{\de}{\delta}
\newcommand{\del}{\partial}

\newcommand{\Om}{\Omega}

\newcommand{{\loc}}{{\mathrm loc}}

\def\meanint{{\diagup\hskip -.42cm\int}}

\begin{document}

\title{Asymptotics for Solutions of Elliptic Equations in Double Divergence Form}
\date{March 13, 2006}
\author{Vladimir Maz'ya}
\address{Link\"oping University, Ohio State University, University of Liverpool}
\author{ Robert McOwen}
\address{Northeastern University}

\begin{abstract}
We consider  weak solutions of an elliptic equation of the form $\del_i\del_i(a_{ij}u)=0$ and their  asymptotic properties at an interior point. We assume that the coefficients are bounded, measurable, complex-valued functions that  stabilize as $x\to 0$ in that the norm of the matrix $(a_{ij}(x)-\de_{ij})$ on the annulus
$B_{2r}\backslash B_r$ is bounded by a function $\Om(r)$, where $\Om^2(r)$ satisfies the Dini condition at $r=0$, as well as some technical monotonicity conditions; under these assumptions, solutions need not be continuous. Our main result is an explicit formula for 
the leading asymptotic term for solutions with at most a mild singularity at $x=0$.
As a consequence, we obtain upper and lower estimates for the $L^p$-norm of solutions, as well
as necessary and sufficient conditions for solutions to be bounded or tend to zero in
$L^p$-mean as $r\to 0$.

\end{abstract}

\maketitle

\newtheorem{theorem}{Theorem}
\newtheorem{lemma}{Lemma}
\newtheorem{proposition}{Proposition}
\newtheorem{corollary}{Corollary}
\newtheorem{definition}{Definition}
\newtheorem{remark}{Remark}
\newtheorem*{main}{Main Theorem}

\addtocounter{section}{-1}
 \section{Introduction}
\smallskip
We are interested in the local behavior of weak solutions to the elliptic equation in
``double divergence form''
\begin{equation}
{\mathcal A}u:=\del_i\del_j(a_{ij}(x)u(x))=0,
\label{eq:pde}
\end{equation}
where we have used  $\del_i=\del/\del x_i$ and the summation convention;  the coefficients $a_{ij}=a_{ji}$ are bounded, measurable, complex-valued functions in a domain to be specified. 
The operator ${\mathcal A}$ arises naturally as 
the formal adjoint ${\mathcal L}^*$ of the operator in ``non-divergent form,'' 
\begin{equation}
{\mathcal L}=\overline{a}_{ij}(x)\del_i\del_j.
\label{eq:L}
\end{equation}
Solutions of (\ref{eq:pde}) are not only important for the solvability of ${\mathcal L} u=f$, but for properties of the Green's function for ${\mathcal L}$. When the coefficients $a_{ij}$ are real-valued functions, the operators ${\mathcal L}$ and ${\mathcal A}$ have been studied by  Sj\"ogren (\cite{Sjo}), Bauman (\cite{Bau1}, \cite{Bau2}, \cite{Bau3}), 
Fabes and Stroock (\cite{FS}), 
Fabes, Garofalo, Mar\'in-Malav\'e, and Salsa (\cite{FGMS}),
 Escauriaza and Kenig (\cite{EK}),
 and Escauriaza (\cite{E1},  \cite{E2}); these papers use techniques, such as the maximum principle, that rely on the coefficients being real-valued. Not only do our techniques apply to complex-valued coefficients, but they give additional information for real-valued coefficients; we shall 
explain this in some detail  at the end of this Introduction. 

We want to study weak solutions of (\ref{eq:pde}) in a neighborhood of an interior point of the domain, say $x=0$, where the coefficients $a_{ij}$ approach $\delta_{ij}$ in the sense that
\begin{equation}
\sup_{r<|x|<2r}|\,{\bf a}(x)-{\bf I}\,|\leq
\Om(r);
\label{eq:aij-asym}
\end{equation}
here ${\bf a}(x)$ is the matrix $(a_{ij})$, ${\bf I}$ is the identity matrix, $|\cdot|$ denotes the matrix norm,  and $\Om(r)\to 0$ as $r\to 0$ in a manner that we shall describe. We remark that, when the coefficients are real-valued, the more general case obtained by replacing $\de_{ij}$ by  constants $\alpha_{ij}$ that satisfy the ellipticity condition can be reduced to (\ref{eq:aij-asym}) by means of an affine change of the $x$ variables. Of course, this reduction of the more general case to (\ref{eq:aij-asym}) is not available when the constants $\alpha_{ij}$ are complex-valued, but we have chosen to treat the special case $\alpha_{ij}=\delta_{ij}$ in order to take advantage of technical simplifications in the formulations and proofs of our results.

The specific hypotheses that we impose on the function $\Om(r)$ in (\ref{eq:aij-asym}) are as follows:
\begin{equation}
\int_0^1 \frac{\Om^2(t)}{t}dt < \infty,
\label{eq:Om1}
\end{equation}
\begin{equation}
\Om(r)\,r^{-1+\e}\quad\hbox{is nonincreasing for $0<r<1$ and},
\label{eq:Om2}
\end{equation}
\begin{equation}
\Om(r)\,r^{n-\e}\quad\hbox{is nondecreasing for $0<r<1$};
\label{eq:Om3}
\end{equation}
here $\e>0$.
Clearly, (\ref{eq:Om1}) together with (\ref{eq:Om2}) or (\ref{eq:Om3}) implies that $\Om(r)\to 0$ as $r\to 0$, so the coefficients $a_{ij}$ are approaching $\de_{ij}$ as $x\to 0$, although perhaps at a slow rate. 

A weak solution of (\ref{eq:pde}) in a domain $U\subset\RR^n$ is a function $u\in L^1_{\loc}(U)$ that satisfies
\begin{equation}
\int_{U} a_{ij}(x)\, u(x)\,\del_j\del_i\eta(x)\,dx=0
\quad\hbox{for all}\ \eta\in C_0^\infty(U).
\label{eq:weaksoln}
\end{equation}
Weak solutions of (\ref{eq:pde}) need not be continuous under our assumptions on the coefficients, so to measure growth or decay as $x\to 0$, we will use the mean in $L^p$ for some $p\in (1,\infty)$:
\begin{equation}
M_p(w,r):=\left(\meanint_{A_r}|w|^p\,dx\right)^{1/p},
\label{eq:Mp}
\end{equation}
where $A_r$ is the annulus $B_{2r}\backslash B_r$ with $B_r=\{x\in \RR: |x|<r\}$;
here (and elsewhere in this paper) the slashed integral denotes the mean value. (We will also use
the notation $M_p(w,r)$ when $w$ is vector or matrix valued; in this case, $|w|$ denotes the norm of $w$.) 

For our local results, we will consider (\ref{eq:pde}) in the unit ball $B_1$ and we will
assume
\begin{equation}
\int_0^1 \frac{\Om^2(t)}{t}dt < \de,
\label{eq:Om4}
\end{equation}
where $\de$ is sufficiently small. In fact, this represents no additional assumption on $\Om(r)$ since we could replace $B_1$ in what follows by a very small ball $B_\gamma$ in order to make the integral  $\int_0^\gamma \Om^2(t)t^{-1}dt$ as small as necessary. 

At times it will be useful to consider solutions of (\ref{eq:pde}) in all of $\RR^n$; in that case,
we assume that $a_{ij}=\de_{ij}$ outside of $B_1$. Our first result concerns such a solution.

\begin{theorem} Let $n\geq 2$, $p\in (1,\infty)$, and $\Om(r)$ satisfy (\ref{eq:Om2}), (\ref{eq:Om3}), and (\ref{eq:Om4}).
There exists a weak solution $Z\in L^p_{\loc}(\RR^n)$ of equation (\ref{eq:pde}) in $\RR^n$ satisfying
\begin{equation}
Z(x)=\exp\left[
-\frac{1}{|\del B_1|}\int_{B_1\backslash B_{|x|}}(a_{ii}(y)-na_{ij}(y)y_iy_j|y|^{-2})\frac{dy}{|y|^n}\right](1+\zeta(x)),
\label{eq:Z-asym}
\end{equation}
where  $\zeta$ satisfies
\[
M_p(\zeta,r)\leq c \, \max\left(\Om(r),\int_0^r\frac{\Om^2(t)}{t}\,dt\right) \quad\hbox{for $0<r<1$.}
\] 
\end{theorem}

\begin{remark}
The proof of Theorem 1 shows that $Z$ has a limit at infinity, $Z(\infty)$, satisfying
\[
|Z(x)-Z(\infty)|\leq c\,\sqrt{\de}\,|x|^{-n}.
\]
But we are more interested in the behavior of $Z$ at the origin; Corollary 3 below can be used to show that $Z$ has at most a mild singularity at the origin (cf.\ (\ref{eq:Mp(u)-est2})).
\end{remark}

Our second theorem uses the solution $Z$ from Theorem 1 to characterize the asymptotics (as $x\to 0$) of weak solutions of (\ref{eq:pde}); because this is a local result, we consider a solution in $B_1$.

\begin{theorem} Let $n>2$, $p\in (1,\infty)$,  and $\Om(r)$ satisfy (\ref{eq:Om2}), (\ref{eq:Om3}), and (\ref{eq:Om4}).
Suppose that $u\in L_{\loc}^p(\overline{B}_1\backslash\{0\})$ is a weak solution of (\ref{eq:pde}) in $B_1$ subject to
the growth condition 
\begin{equation}
M_p(u,r)\leq c\,r^{2-n+\e_0},
\label{eq:u-bd}
\end{equation}
where $\e_0>0$. Then there exists a constant $C$ (depending on $u$) such that
\begin{equation}
u(x)=CZ(x)+w(x),
\label{eq:u-asym}
\end{equation}
where the remainder term $w$ satisfies
\begin{equation}
M_p(w,r)\leq c\,r^{1-\e_1}
\label{eq:w-bd}
\end{equation}
for $0<r<1$ and any $\e_1>0$.
\end{theorem}

\begin{remark} The restriction $n>2$ in Theorem 2 is caused by the existence of solutions for the Laplacian with logarithmic growth at $x=0$ when $n=2$. A refinement of the techniques used in proving Theorem 2 would be required to cover the case $n=2$.
\end{remark}

The following two results are immediate consequences of Theorems 1 and 2.

\begin{corollary}
Under the hypotheses of Theorem 2, the condition
\begin{equation}
\liminf_{r\to 0} \int_{B_1\backslash B_{r}}{\rm Re}(a_{ii}(y)-na_{ij}(y)y_iy_j|y|^{-2})\frac{dy}{|y|^n}>-\infty
\end{equation}
is necessary and sufficient for the boundedness of $M_p(u,r)$ as $r\to 0$ for all
$u\in L_{\loc}^p(\overline{B}_1\backslash\{0\})$ which are weak solutions
 of (\ref{eq:pde}) and satisfy (\ref{eq:u-bd}).
\end{corollary}

\begin{corollary}
Under the hypotheses of Theorem 2, the condition
\begin{equation}
\lim_{r\to 0} \int_{B_1\backslash B_{r}}{\rm Re}(a_{ii}(y)-na_{ij}(y)y_iy_j|y|^{-2})\frac{dy}{|y|^n}=+\infty
\end{equation}
is necessary and sufficient for  $M_p(u,r)\to 0$ as $r\to 0$ for all
$u\in L_{\loc}^p(\overline{B}_1\backslash\{0\})$ which are weak solutions
 of (\ref{eq:pde}) and satisfy (\ref{eq:u-bd}).
\end{corollary}

\noindent
Moreover, we can use Theorems 1 and 2 to derive upper and lower estimates for the
$L^p$-mean of solutions.
Since $({\rm Re}(a_{ij})-\de_{ij})$ is a symmetric real matrix, it is not difficult to verify that
\begin{equation}
-2(n-1)|\,{\rm Re}\,{\bf a}(y)-{\bf I}\,|\leq {\rm Re}(a_{ii}(y)-na_{ij}(y)y_iy_j|y|^{-2})\leq 
2(n-1)|\,{\rm Re}\,{\bf a}(y)-{\bf I}\,|.
\end{equation} 
Using $|\,{\rm Re}\,{\bf a}(y)-{\bf I}\,|\leq |\,{\bf a}(y)-{\bf I}\,|$ and (\ref{eq:aij-asym}), we find that
(\ref{eq:Z-asym}) also yields the following upper and lower estimates.

\begin{corollary}
The solution $Z$ in Theorem 1 satisfies for $ r\in (0,1)$
\begin{equation}
c_1\exp\left(-2(n-1)\int_r^1\Om(s)\frac{ds}{s}\right)\leq M_p(Z,r)\leq c_2\exp\left(2(n-1)\int_r^1\Om(s)\frac{ds}{s}\right),
\label{eq:Mp(u)-est1}
\end{equation} 
where $c_1$, $c_2$ are positive constants.
\end{corollary}

\noindent
Even when the coefficients $a_{ij}$ are real, the upper and lower bounds (\ref{eq:Mp(u)-est1}) appear to be new.  

The principal analytic content of our results is contained in Theorem 1. The method of its proof is independent of, but related to, the asymptotic theory developed in ~\cite{KM3}. In particular, $L_p$-means of type (\ref{eq:Mp}) were extensively used in ~\cite{KM1} and ~\cite{KM2}. The asymptotic formula that we obtain is analogous to that of ~\cite{KM4}, where an asymptotic representation near the boundary was obtained for solutions to the Dirichlet problem for elliptic equations in divergence form with discontinuous coefficients.

Now let us turn to the comparison of our results with the extensive work of the authors cited in the first paragraph; we refer to the excellent exposition in Escauriaza  \cite{E2} for a more detailed description of these results and references to the literature. When ${\mathcal L}$ is uniformly elliptic with real-valued, measurable, and bounded (although not necessarily continuous) coefficients on $\RR^n$, the previous work shows: i) any nonnegative weak solution $Z\in L^1_{\loc}(B_1)$ of (\ref{eq:pde}) in $B_1$ satisfies  $Z\in L^{n/(n-1)}(B_1)$ (cf.\ \cite{FS}); ii) for a fixed nontrivial nonnegative weak solution $Z\in L^1_{\loc}(B_1)$ of (\ref{eq:pde}) in $B_1$, and every weak solution $u\in L^1_{\loc}(B_1)$ of (\ref{eq:pde}) in $B_1$, the function $u/Z$ is H\"older continuous in $B_1$ (cf.\ \cite{Bau1}); and iii) the existence of a unique nonnegative weak solution $Z$ of (\ref{eq:pde}) in $\RR^n$ that satisfies $\int_{B_1} Z\,dx=|B_1|$ and may be used to estimate the Green's function for the operator ${\mathcal L}$ on $\RR^n$ (cf.\ \cite{E2}).  These previous results are quite general, but do not apply to the case of complex coefficients that we consider because they rely upon the maximum principle. Moreover, even when the coefficients are real-valued, our results are somewhat different in nature than the previous ones: instead of estimates, we have obtained 
an asymptotic description of $Z$ near a point where the $a_{ij}$ are continuous in the sense of  (\ref{eq:aij-asym}). For a more direct comparison, when the $a_{ij}$ are real-valued and continuous at $0$, Escauriaza  \cite{E2} has obtained upper and lower estimates for the $L^1$-norm of $Z$: for any $\gamma>0$ there exists a constant $N_\gamma$ such that
\begin{equation}
N_\gamma^{-1} r^{\gamma}\leq \meanint_{B_r} Z\,dx \leq N_\gamma r^{-\gamma} \quad\hbox{for}\ 0<r<1.
\label{eq:EscauriazaEst}
\end{equation}
When the $a_{ij}$ satisfy our stronger sense of continuity (\ref{eq:aij-asym}), the explicit upper and lower  bounds  (\ref{eq:Mp(u)-est1}) imply $L^p$-bounds analogous to (\ref{eq:EscauriazaEst}).
In fact, as we shall see in Section 2, the assumptions on $\Omega(r)$ imply that $\Om(r)\leq c\sqrt{\de}$ for $0<r<1$ (cf.\ (\ref{eq:Om^2<de})), so we easily obtain from  (\ref{eq:Mp(u)-est1})
\begin{equation}
N_\gamma^{-1} r^{\gamma}\leq M_p(Z,r)\leq N_\gamma r^{-\gamma} \qquad\hbox{for} \ 0<r<1.
\label{eq:Mp(u)-est2}
\end{equation} 
More importantly, however, our necessary and sufficient conditions for solutions to be  bounded in $L^p$-mean as $r\to 0$ (Corollary 1) or tend to zero in 
 $L^p$-mean as $r\to 0$ (Corollary 2) have not been obtained previously.

\section{Preliminary Estimates}

In this and the next section, we will use the spherical mean of a function $w$. For notational
convenience, we denote the spherical mean using an ``overbar'':
\begin{equation}
\overline{w}(r)=\meanint_{\del B_1} w(r\theta)\,ds.
\label{eq:sph-mean}
\end{equation}
This should cause no confusion with complex conjugation since we will not have occasion to use the latter in these two sections.
In particular, in this section we are concerned with solving an equation of the form
\begin{equation}
-\Lap v=\del_i\del_j(F_{ij})-\overline{\del_i\del_j(F_{ij})} \quad\hbox{in}\ \RR^n.
\label{eq:Lap}
\end{equation}
Here $F_{ij}\in L^1_{\loc}(\RR^n)$ and derivatives are interpreted in the sense of distributions.
The norm of the matrix ${\mathcal F}=(F_{ij})$ will be denoted by $|{\mathcal F}|$.

\begin{proposition}
Suppose that $F_{ij}\in L^p_{\loc}(\RR^n\backslash\{0\})$ satisfies 
\begin{equation}
\int_{|x|<1}|{\mathcal F}(x)|\,dx+\int_{|x|>1}|{\mathcal F}(x)||x|^{-n-1}dx<\infty.
\label{eq:F-cond}
\end{equation}
Then there exists a weak solution $v\in L^p_{\loc}(\RR^n\backslash\{0\})$ of (\ref{eq:Lap})
that satisfies
\begin{equation}
M_p(v,r)\leq 
c\left(\tilde{M}_p({\mathcal F},r)+r\int_{|x|>r}|{\mathcal F}(x)||x|^{-n-1}dx
+r^{-n}\int_{|x|<r}|{\mathcal F}(x)|dx\right),
\label{eq:Lap-est1}
\end{equation}
where $c$ is independent of $r$ and we have introduced
\[
\tilde{M}_p(w,r):=\left(\meanint_{r/2<|x|<4r}|w(x)|^p\,dx\right)^{1/p}.
\]
\label{pr:Lap}
\end{proposition}
\noindent
{\bf Proof:}
It suffices to prove the result for $F_{ij}\in C_0^\infty(\RR^n\backslash\{0\})$ since the general case can be handled by a standard approximation argument.
The function $v$ is defined by convolution with $\Gamma$, the fundamental solution for the Laplacian:
\[
v=\Gamma\star(\del_i\del_j F_{ij}-\overline{\del_i\del_j F_{ij}}).
\]
Using $\int_{\RR^n} f(y)\overline g(|y|)\, dy=\int_{\bf R^n} \overline f(|y|) g(y)\,dy$, we can write this as
\[
v(x)=\int_{\RR^n}\left(\Gamma(|x-y|)-\overline{\Gamma(|x-\cdot|)}(|y|)\right)\del_i\del_j F_{ij}(y)\,dy.
\]
Now to compute the spherical mean of the fundamental solution, we can use the mean value theorem for harmonic functions to conclude that 
\[
\overline{\Gamma(|x-\cdot|)}(|y|) = \Gamma(\max\{|x|,|y|\}).
\]
This enables us to express $v$ as
\[
v(x) =\frac{\del^2}{\del x_i\del x_j}\int_{\RR^n}\Gamma(|x-y|)F_{ij}(y)\,dy-
\Gamma(|x|)\int_{|y|<|x|}\frac{\del^2 F_{ij}}{\del y_i\del y_j}\,dy-
\int_{|y|>|x|}\Gamma(|y|)\frac{\del^2 F_{ij}}{\del y_i\del y_j}\,dy.
\]
Now, integration by parts yields
\[
\int_{|y|<|x|}\frac{\del^2 F_{ij}}{\del y_i\del y_j}\,dy=\int_{|y|=|x|}\frac{\del F_{ij}}{\del y_j}\frac{y_i}{|y|}\,dS_y,
\]
and
\[
\int_{|y|>|x|}\Gamma(|y|)\frac{\del^2 F_{ij}}{\del y_i\del y_j}\,dy
=-\int_{|y|>|x|}\Gamma'(|y|)\frac{y_i}{|y|}\frac{\del F_{ij}}{\del y_j}\,dy-
\int_{|y|=|x|}\Gamma(|y|)\frac{\del F_{ij}}{\del y_j}\frac{y_i}{|y|}\,dS_y
\]
\[
=\int_{|y|>|x|}\frac{\del}{\del y_j}\left(\Gamma'(|y|)\frac{y_i}{|y|}\right)F_{ij}(y)\,dy
+\int_{|y|=|x|}\!\left(\Gamma'(|y|)\frac{y_j}{|y|}F_{ij}(y)-\Gamma(|y|)\frac{\del F_{ij}}{\del y_j}\right)\!
\frac{y_i}{|y|}\,dS_y.
\]
Consequently,
\[
v(x)= \ \frac{\del^2}{\del x_i\del x_j}\int_{\RR^n}\Gamma(|x-y|)F_{ij}(y)\,dy
-\int_{|y|>|x|}\left(\frac{\del^2}{\del y_j\del y_i}\Gamma(|y|)\right)F_{ij}(y)\,dy
\]
\[
\ -\Gamma'(|x|)\int_{|y|=|x|}\frac{y_iy_j}{|y|^2}F_{ij}(y)\,dS_y.
\]

Now introduce $\chi_0$ and $\chi_\infty$ as the characteristic functions of $B_{r/2}$ 
 and $B_{4r}^c$, and let $\chi_1=1-\chi_0-\chi_\infty$ be the characteristic function of the annulus 
 $\tilde A_r:=B_{4r}\backslash B_{r/2}$. Then
 \[
 v(x)-\frac{\del^2}{\del x_i\del x_j}\int_{\RR^n}\Gamma(|x-y|)(\chi_1F_{ij})(y)\,dy=
 \]
 \[
\int_{\RR^n}\frac{\del^2}{\del x_i\del x_j}\Gamma(|x-y|)(\chi_0F_{ij})(y)\,dy+
\int_{\RR^n}\frac{\del^2}{\del y_i\del y_j}(\Gamma(|x-y|)-\Gamma(|y|))(\chi_\infty F_{ij})(y)\,dy
\]
\[
-\int_{B_{4r}\backslash B_{|x|}}\left(\frac{\del^2}{\del y_i\del y_j}\Gamma(|y|)\right)F_{ij}(y)\,dy
-\Gamma'(|x|)\int_{|y|=|x|}\frac{y_iy_j}{|y|^2}F_{ij}(y)\,dS_y.
\]
We can estimate the four integral kernels and obtain that the right hand side is bounded by
\[
c\left(\frac{1}{|x|^n}\int_{B_{r/2}}|F_{ij}(y)|dy+|x|\int_{B_{4r}^c}\frac{|F_{ij}(y)|dy}{|y|^{n+1}}
+|x|\int_{B_{4r}\backslash B_{|x|}}\frac{|F_{ij}(y)|dy}{|y|^{n+1}}+\overline{|F_{ij}|}(|x|)
\right)
\]
This provides us with the following pointwise bound:
\[
\left|v(x)- \frac{\del^2}{\del x_i\del x_j} \int_{\RR^n}\Gamma(|x-y|)(\chi_1F_{ij})(y)\,dy\right|  
\]
\begin{equation}
\leq c  \left(\overline{|F_{ij}|}(|x|)+  |x|\int_{B_r^c}|F_{ij}(y)|\frac{dy}{|y|^{n+1}}+|x|^{-n}\int_{B_r}|F_{ij}(y)|dy\right).
\label{eq:ptwse-bd}
\end{equation}

Using the $L^p$-boundedness of singular integral operators on $\RR^n$ (see ~\cite{St}), we have
\[
\left\|\frac{\partial^2}{\partial x_i\partial x_j}
\int_{{\bf R}^n}\Gamma(|x-y|)(\chi_1 F_{ij})(y)dy\right\|_{L^p(A_r)}  \leq 
 \ \|\chi_1 {\mathcal F}\|_{L^p({\bf R}^n)}=\|{\mathcal F}\|_{L^p(\tilde{A}_r)}.
\]
Elementary estimates  may be applied to the remaining terms in (\ref{eq:ptwse-bd}) to obtain (\ref{eq:Lap-est1}), completing the proof.
\hfill $\Box$

\medskip
The integrals in (\ref{eq:Lap-est1}) can be estimated in terms of $M_p$. For example, we substitute
\[
|x|^{-n-1}=c_n\int_{\frac{|x|}{2}<|y|<|x|}|y|^{-2n-1}\,dy,
\]
into the first integral and change order of integration to obtain the estimate 
\[
r\int_{|x|>r}|{\mathcal F}(x)||x|^{-n-1}dx\leq 
c \,r\int_{|y|>r/2}\int_{|y|<|x|<2|y|}|{\mathcal F}(x)|\,dx\, \frac{dy}{|y|^{2n+1}}
\leq c\,r\int_{r/2}^\infty M_p({\mathcal F},\rho)\,\frac{d\rho}{\rho^2}.
\]
Similarly, we can show
\[
r^{-n}\int_{|x|<r}|{\mathcal F}(x)|\,dx\leq c\, r^{-n}\int_0^r M_p({\mathcal F},\rho)\,\rho^{n-1}\,d\rho
\]
and
\[
\tilde M_p({\mathcal F},r)^p\leq c r^{-n}\int_{r/2}^{4r} M_p({\mathcal F},\rho)^p\,\rho^{n-1}\,d\rho.
\]
Elementary estimates show that terms involving integration over $r/2<\rho<r$ and $2r<\rho<4r$ can be respectively dominated  by the terms involving integration over $0<\rho<r$ and $\rho>r$, so we obtain the following.

\begin{corollary}
Under the hypotheses of Proposition \ref{pr:Lap}, the weak solution $v$ obtained there satisfies
\begin{equation}
M_p(v,r)\leq c\,\left(r\int_r^\infty M_p({\mathcal F},\rho)\rho^{-2}\,d\rho+
r^{-n}\int_0^rM_p({\mathcal F},\rho)\rho^{n-1}\,d\rho\right).
\label{eq:Lap-est2}
\end{equation}
\label{corollary}
\end{corollary}
\noindent

\section{Proof of Theorem 1}

We shall prove Theorem 1 by reducing the problem of finding $Z$ to solving an operator equation of the form $(I+\tilde T)V=f$, where $V$ and $f$ are elements of a Banach space $X$ of functions on $\RR^n\backslash\{0\})$, and $\tilde T$ is an integral operator of small norm on $X$. However, this reduction will take a few steps. To begin, let $r=|x|$, $\theta=x/|x|$, and $\eta\in C_0^\infty((0,\infty))$ be arbitrary. For $Z\in L_{\loc}^p(\RR^n)$ to be  a weak solution of (\ref{eq:pde}), we must have
\[
0=\int_{\RR^n}\del_i\del_j\eta(|x|)a_{ij}(x)Z(x)\,dx
\]
\[
=\int_0^\infty \left(\eta''(r)\int_{\del B_1} Z(r\theta)a_{ij}(r\theta)\theta_i\theta_j ds+
\frac{\eta'(r)}{r}\int_{\del B_1} Z(r\theta)(a_{ii}(r\theta)-a_{ij}(r\theta)\theta_i\theta_j)\,ds\right)r^{n-1}dr,
\]
where $ds$ denotes surface measure on the unit sphere, $\del B_1$.
Hence,
\[
0=\int_0^\infty \eta'(r)\left(-\frac{d}{dr}\left[r^{n-1}\int_{\del B_1}Z(r\theta)a_{ij}(r\theta)\theta_i\theta_j\,ds\right] \right.
\]
\[
\left. +r^{n-2}\int_{\del B_1} Z(r\theta)(a_{ii}(r\theta)-a_{ij}(r\theta)\theta_i\theta_j)ds\right)dr,
\]
where the derivative is understood in the distributional sense. This implies
\begin{equation}
-r^{n-1}\frac{d}{dr}\int_{\del B_1} Z(r\theta)a_{ij}(r\theta)\theta_i\theta_jds
+r^{n-2}\int_{\del B_1} Z(r\theta) (a_{ii}(r\theta)-n a_{ij}(r\theta)\theta_i\theta_j)ds=C,
\label{eq:Z-eq}
\end{equation}
where $C$ is an arbitrary constant. In what follows, we will take $C=0$; as we shall see, the solution that we construct will in fact be a weak solution of (\ref{eq:pde}) on all of $\RR^n$, not just $\RR^n\backslash\{0\}$. (See also the Remark at the end of this section.)

Let us introduce
\begin{equation}
v(r\theta):= Z(r\theta)-\overline{Z}(r),
\label{eq:def-v}
\end{equation}
where $\overline{Z}$ is the spherical mean as in (\ref{eq:sph-mean}).
We may now express (\ref{eq:Z-eq}) as
\begin{equation}
y'(r)+\frac{Q(r)}{r}y(r)=\frac{1}{r}Kv(r),
\label{eq:ode-y}
\end{equation}
where
\begin{equation}
y(r):=\meanint_{\del B_1} Z(r\theta) a_{ij}(r\theta)\theta_i\theta_j\, ds,
\label{eq:def-y}
\end{equation}
and
\begin{equation}
Q(r):=n-\frac{\alpha_0(r)}{\alpha(r)},
\label{eq:def-Q}
\end{equation}
with
\begin{equation}
\alpha_0(r):=\meanint_{\del B_1} a_{ii}(r\theta)\, ds \quad\hbox{ and}\quad
\alpha(r):=\meanint_{\del B_1} a_{ij}(r\theta)\theta_i\theta_j\,ds;
\label{eq:alphas}
\end{equation}
in (\ref{eq:ode-y}) we also have used
\begin{equation}
Kv(r):=\meanint_{\del B_1} v(r\theta) a_{ii}(r\theta)\,ds-
\frac{\alpha_0(r)}{\alpha(r)}
\meanint_{\del B_1} v(r\theta)a_{ij}(r\theta)\theta_i\theta_j\, ds.
\label{eq:def-K}
\end{equation}
It follows from (\ref{eq:aij-asym})  that $|\alpha_0(r)-n|\leq c\,\Om(r)$, $|\alpha(r)-1|\leq c\,\Om(r)$, and 
\begin{equation}
|Q(r)|\leq c\,\Om(r),
\label{eq:Q-Om}
\end{equation}
so $Q(r)\to 0$ as $r\to 0$.
Since $\overline{v}(r)=0$, we can also write (\ref{eq:def-Q}) as
\[
Kv(r)=
\meanint_{\del B_1}  v(r\theta)(a_{ii}(r\theta)-n)\,ds
-\frac{\alpha_0(r)}{\alpha(r)}
\meanint_{\del B_1}v(r\theta)(a_{ij}(r\theta)-\de_{ij})\theta_i\theta_j\,ds.
\]
In this last form it is evident that $K$ satisfies
\begin{equation}
M_p(Kv,r)\leq c\,\Om(r) \, M_p(v,r).
\label{eq:Kv-bd}
\end{equation}

To obtain another equation involving $y$ and $v$, we start from the identity
\begin{equation}
\Lap v=\overline{\del_i\del_j((a_{ij}-\de_{ij})v)}-\del_i\del_j((a_{ij}-\de_{ij})v)
+\overline{\del_i\del_j((a_{ij}-\de_{ij})\overline{Z})}- \del_i\del_j((a_{ij}-\de_{ij})\overline{Z}).
\label{eq:pde-v}
\end{equation}
Noting that 
\[
y(r)=\alpha(r)\overline{Z}(r)
+\meanint_{\del B_1} v(r\theta)(a_{ij}(r\theta)\theta_i\theta_j -1)\,ds,
\]
we can rewrite (\ref{eq:pde-v}) as
\begin{equation}
\Lap v=\overline{\del_i\del_j(B_{ij}(v))}-\del_i\del_j(B_{ij}(v))
+ \overline{\del_i\del_j(\phi_{ij}y)}-\del_i\del_j(\phi_{ij}y),
\label{eq:pde2-v}
\end{equation}
where 
\begin{equation}
B_{ij}(v)(x):=(a_{ij}(x)-\de_{ij})\left(v(x)-\frac{1}{\alpha(r)}\meanint_{\del B_1}v(r\theta)(a_{ij}(r\theta)\theta_i\theta_j-1)\,ds
\right)
\label{eq:def-B}
\end{equation}
and
\begin{equation}
\phi_{ij}(x)=\frac{a_{ij}(x)-\de_{ij}}{\alpha(r)}.
\label{eq:def-phi}
\end{equation}
Using (\ref{eq:aij-asym}), it is clear that
\begin{equation}
M_p(B_{ij}(v),r)\leq c\,\Om (r)\, M_p(v,r)
\label{eq:Mp-B}
\end{equation}
and
\begin{equation}
|\phi_{ij}(x)|\leq c\,\Om (r).
\label{eq:phi-bd}
\end{equation}
These estimates indeed hold for $0<r<\infty$, where we extend $\Om$ to be zero for $r>1$.
Inverting the Laplacian by means of the fundamental solution, (\ref{eq:pde2-v}) becomes 
\begin{equation}
v+Sv+Ty=0,
\label{eq:int-eqn1}
\end{equation}
where we may use Corollary \ref{corollary} with (\ref{eq:Mp-B}) to obtain 
\begin{equation}
M_p(Sv,r)\leq c\left(r\int_r^\infty \Om(\rho)M_p(v,\rho)\rho^{-2}\,d\rho
+r^{-n}\int_0^r\Om(\rho) M_p(v,\rho)\rho^{n-1}\,d\rho\right)
\label{eq:Mp-S}
\end{equation}
and with (\ref{eq:phi-bd}) to obtain
\begin{equation}
M_p(Ty,r)\leq c\left(r\int_r^\infty \Om(\rho)M_p(y,\rho)\rho^{-2}\,d\rho
+r^{-n}\int_0^r\Om(\rho) M_p(y,\rho)\rho^{n-1}\,d\rho\right).
\label{eq:Mp-T}
\end{equation}

To simplify the equations, let us
introduce the function
\begin{equation}
E(r)=\exp\left(\int_r^\infty \frac{Q(t)}{t}\,dt\right),
\label{eq:E}
\end{equation}
where there is no problem with convergence of the integral since $Q(t)=0$ for $t>1$.
Notice that $E(r)$ is continuous for $r\in (0,\infty)$, $E(r)\equiv 1$ for $r\geq 1$, and for any $r,\rho\in (0,\infty)$ we have
\begin{equation}
E^{-1}(r)E(\rho)=\exp\left(\int_\rho^r \frac{Q(t)}{t}\,dt\right).
\label{eq:E-1*E}
\end{equation}
To derive an estimate for this expression, let us use (\ref{eq:Om2}) with $\e\in (0,1/2)$ and (\ref{eq:Om4})   to conclude
\begin{equation}
\de\geq\int_{r/2}^r \frac{\Om^2(t)}{t}dt\geq \Om^2(r)r^{-2+2\e}\int_{r/2}^r t^{1-2\e}dt\geq c\,\Om^2(r).
\label{eq:Om^2<de}
\end{equation}
As a consequence of (\ref{eq:Q-Om}), we have
\begin{equation}
\left(\frac{\rho}{r}\right)^{c\sqrt{\de}}\leq 
\left|\,\exp\left(\pm\int_\rho^r \frac{Q(t)}{t}\,dt\right)\right|
\leq \left(\frac{r}{\rho}\right)^{c\sqrt{\de}}
\quad\hbox{for}\ 0<\rho\leq r\leq 1.
\label{eq:intQ-est}
\end{equation}

Now let us express equations (\ref{eq:ode-y}) and (\ref{eq:int-eqn1}) in terms of the new dependent variables
\begin{equation}
Y(r)=E^{-1}(r)y(r)
\quad\hbox{and}\quad
V(x)=E^{-1}(|x|)v(x).
\label{eq:def-Y,V}
\end{equation}
Since the operator $K$ only involves integration in $\theta$, it is clear that
\begin{equation}
Kv(r)=E(r)KV(r),
\label{eq:EKV}
\end{equation}
and so the equation (\ref{eq:ode-y}) can be expressed as
\begin{equation}
Y'(r)=\frac{1}{r} KV(r).
\label{eq:ode-Y}
\end{equation}
We will be assuming below that $M_p(V,r)\leq c\,\Omega(r)$ for $0<r<1$, so 
$M_p(KV,r)\leq c\,\Om^2(r)$, enabling us to
 integrate (\ref{eq:ode-Y}) to obtain
\begin{equation}
Y(r)=Y(0)+\int_0^r\frac{KV(\rho)}{\rho}\,d\rho.
\label{eq:sol-Y}
\end{equation}
The equation (\ref{eq:int-eqn1}), on the other hand, is replaced by
\begin{equation}
V+{S}_1V+{T}_1Y=0,
\label{eq:int-eqn2}
\end{equation}
where
\begin{equation}
{S}_1:=E^{-1}SE
\quad\hbox{and}\quad
{T}_1:=E^{-1}TE,
\label{eq:def-ST}
\end{equation}
with $E$ representing the multiplication operator defined by the function (\ref{eq:E}).

Now let us substitute (\ref{eq:sol-Y}) into (\ref{eq:int-eqn2}) to finally obtain the operator equation that we want to solve:
\begin{equation}
V+{S}_1V+{T}_2V=-Y(0){T}_1(1),
\label{eq:int-eqnV}
\end{equation}
where 
\begin{equation}
T_2V(x)=E^{-1}(|x|)T_{\rho\to x}\left[E(\rho)\int_0^\rho KV(t)\,dt\right].
\label{eq:T2}
\end{equation}
(In (\ref{eq:int-eqnV}), notice that ${T}_1$ operates on functions of a single variable, say $\rho$, and 
${T}_1(1)$ represents the action of ${T}_1$ on the function that is identically $1$.) For a given choice of $Y(0)$, we can solve (\ref{eq:int-eqnV}) uniquely for $V$ provided we can show that the two integral operators involved have small norm on an appropriate function space. Consider the functions $w$ in $L^p_{\loc}(\RR\backslash\{0\})$ for which the norm
\begin{equation}
\|w\|_{p,\Om}:=\sup_{0<r<1}\frac{M_p(w,r)}{\Om(r)}+\sup_{r>1}\frac{M_p(w,r)}{\sqrt{\de}\,r^{-n}}
\label{eq:def-X}
\end{equation}
is finite, and take the closure to form a Banach space $X$.
We want to show that the right hand side of (\ref{eq:int-eqnV}) is in $X$ and that the integral operators $S_1$ and $T_2$ map $X$ to itself with small norm. It will be useful to observe that the continuity and positivity of $|E(r)|$ implies that, for any $w\in L^p_{\loc}(\RR\backslash\{0\})$, we have
\begin{equation}
M_p(Ew,r)=|E(\tilde r)| M_p(w,r) \quad\hbox{for some}\ \tilde r=\tilde r_w\in (r,2r),
\label{eq:Mp-E}
\end{equation}
with an analogous statement for $E^{-1}$.

To show $T_1(1)\in X$, we must estimate $M_p(T_1(1),r)$ separately for $0<r<1$ and $r>1$.
For $0<r<1$, we can use (\ref{eq:Mp-E}) to find $\tilde r\in (r,2r)$ so that
\[
M_p(T_1(1),r)=|E^{-1}(\tilde r)| M_p(T_{\rho\to r}[E(\rho)],r),
\]
and then use (\ref{eq:Mp-T}) to estimate
\[
M_p(T_1(1),r)\leq c \, |E^{-1}(\tilde r)|\left( r\int_r^\infty\Om(\rho)M_p(E,\rho)\rho^{-2}\,d\rho+
r^{-n}\int_0^r\Om(\rho)M_p(E,\rho)\rho^{n-1}\,d\rho\right)
\]
\[
=c \, |E^{-1}(\tilde r)|\left( r\int_r^\infty\Om(\rho)|E(\tilde\rho)|\rho^{-2}\,d\rho+
r^{-n}\int_0^r\Om(\rho)|E(\tilde\rho)|\rho^{n-1}\,d\rho\right),
\]
where $\tilde\rho\in(\rho,2\rho)$ by (\ref{eq:Mp-E}).
But now we can use (\ref{eq:E-1*E}) and (\ref{eq:intQ-est}) to conclude
\[
M_p(T_1(1),r)\leq c\left(r\int_r^{1}\Om(\rho)(\rho/r)^{c\sqrt{\de}}\rho^{-2}\,d\rho
+r^{-n}\int_0^r\Om(\rho)(r/\rho)^{c\sqrt{\de}}\rho^{n-1}\,d\rho
\right)
\]
Using the monotonicty properties (\ref{eq:Om2}) and (\ref{eq:Om3}), we obtain
\[
M_p(T_1(1),r)\leq c\left(\Om(r)r^{\e-c\sqrt{\de}}\int_r^1\rho^{-\e-1+c\sqrt{\de}}d\rho
+\Om(r)r^{c\sqrt{\de}-\e}\int_0^r\rho^{\e-c\sqrt{\de}-1}d\rho\right).
\]
Provided $\de$ is sufficiently small that $\e-c\sqrt{\de}>0$, we conclude that
\begin{equation}
M_p(T_1(1),r)\leq c\,\Om(r)\quad\hbox{for}\ 0<r<1.
\label{eq:RHS-r<1}
\end{equation}
On the other hand, for $r>1$ we use $E^{-1}(r)\equiv 1$ and $\Om(r)\equiv 0$ to estimate
\[
M_p(T_1(1),r)= M_p(T_{\rho\to r}[E(\rho)],r)\leq c \,r^{-n}\int_0^1\Om(\rho)M_p(E,\rho)\rho^{n-1}\,d\rho
\]
\[
=c \,r^{-n}\int_0^1\Om(\rho)|E(\tilde\rho)|\rho^{n-1}\,d\rho
\leq c\,r^{-n}\int_0^1 \Om(\rho)\rho^{n-c\sqrt{\de}-1}\,d\rho.
\]
Provided $n-c\sqrt{\de}>0$, we obtain
\begin{equation}
M_p(T_1(1),r)\leq c\,\sqrt{\de}\,r^{-n}\quad\hbox{for}\ r>1.
\label{eq:RHS-r>1}
\end{equation}
The two estimates (\ref{eq:RHS-r<1}) and (\ref{eq:RHS-r>1}) together confirm that $T_1(1)\in X$.

Now let us show that $S_1$ maps $X$ to itself with small operator norm. We suppose that $\|V\|_{p,\Om}\leq 1$ and estimate $M_p(S_1V,r)$ separately for $0<r<1$ and $r>1$. For
$0<r<1$ we have $M_p(V,r)\leq \Om(r)$ and we  can argue as in the previous paragraph to obtain
\[
M_p(S_1V,r)\leq c\left(r\int_r^1 \Om^2(\rho)(\rho/r)^{c\sqrt{\de}}\rho^{-2}\,d\rho+r^{-n}\int_0^r\Om^2(\rho)(r/\rho)^{c\sqrt{\de}}\rho^{n-1}\,d\rho\right).
\]
Using (\ref{eq:Om^2<de}) and the monotonicity analysis as above, we conclude that
\begin{equation}
\frac{M_p(S_1V,r)}{\Om(r)}\leq c\sqrt{\de}\quad\hbox{for}\ 0<r<1.
\label{eq:S1-r<1}
\end{equation}
For $r>1$ we use (\ref{eq:intQ-est}) with $r=1$ and (\ref{eq:Om^2<de}) to obtain
\[
M_p(S_1V,r)\leq c\,r^{-n}\int_0^1\Om^2(\rho)\rho^{n-c\sqrt{\de}-1}\,d\rho
\leq c\,r^{-n}\de\int_0^1\rho^{n-c\sqrt{\de}-1}\,d\rho\leq c\,r^{-n}\de
\]
provided $\de$ is sufficiently small, and we conclude that
\begin{equation}
\frac{M_p(S_1V,r)}{\sqrt{\de}\,r^{-n}}\leq c\sqrt{\de}\quad\hbox{for}\ r>1.
\label{eq:S1-r>1}
\end{equation}
The estimates (\ref{eq:S1-r<1}) and (\ref{eq:S1-r>1}) together show that $S_1$ maps $X$ to itself with small operator norm.

Finally, we estimate $T_2$. For $0<r<1$ and $M_p(V,r)\leq \Om(r)$, we argue as before to write 
\[
M_p(T_2V,r)\leq c \left( r\int_r^\infty \Om(\rho) \, (\rho/r)^{c\sqrt{\de}}\, M_p\left[\int_0^\rho KV(t)\frac{dt}{t},\rho\right]\rho^{-2}\,d\rho\right.
\]
\[
\left.
+\, r^{-n}\int_0^r\Om(\rho)\,(r/\rho)^{c\sqrt{\de}}\,M_p\left[\int_0^\rho KV(t)\frac{dt}{t},\rho\right]\rho^{n-1}\,d\rho\right).
\]
But we can use  (\ref{eq:Kv-bd}) to obtain $M_p(KV,\rho)\leq c\,\Om(\rho)M_p(V,\rho)\leq c\,\Om^2(\rho)$ so that
\[
M_p(T_2V,r)\leq c \,\left( r\int_r^1 \Om(\rho) \, \left(\rho/r\right)^{c\sqrt{\de}}\left[\int_0^\rho\frac{\Om^2(t)}{t}\,dt\right]\,\rho^{-2}\,d\rho\,+\right.
\]
\[
\left. r^{-n}\int_0^r \Om(\rho) \, \left(r/\rho\right)^{c\sqrt{\de}}\left[\int_0^\rho\frac{\Om^2(t)}{t}\,dt\right]\,\rho^{n-1}\,d\rho
\right)
\leq c\,\de\,\Om(r).
\]
Using (\ref{eq:Om^2<de}) and the monotonicity argument, we have 
\begin{equation}
\frac{M_p(T_2V,r)}{\Om(r)}\leq c\,\de\quad\hbox{for}\ 0<r<1.
\label{eq:T2-r<1}
\end{equation}
For $r>1$, we use (\ref{eq:intQ-est}) with $r=1$ and (\ref{eq:Q-Om}) to obtain
\[
M_p(T_2V,r)\leq c\,r^{-n}\int_0^1\Om^2(\rho)\rho^{-c\sqrt{\de}+n-2}\,d\rho
\leq c\,r^{-n}\de\int_0^1\rho^{-c\sqrt{\de}+n-2}\,d\rho\leq c\,r^{-n}\de,
\]
provided $\de$ is sufficiently small.
These estimates show that $T_2$ maps $X$ to itself with small operator norm.

Since  $S_1$ and $T_2$ have small operator norms on $X$, we conclude that 
(\ref{eq:int-eqnV}) admits a unique solution $V\in X$. 
Let us now investigate the implications for the weak solution $Z(x)=\overline{Z}(|x|)+v(x)$ that we are trying to construct. Tracing back through the definitions, we see that our solution of (\ref{eq:pde}) is given by
\begin{equation}
Z(x)=\frac{E(r)}{\alpha(r)}\left(Y(0)+Y_1(r)-\meanint_{\del B_1} V(r\theta)a_{ij}(r\theta)\theta_i\theta_j\,ds+\alpha(r)V(x)\right),
\label{eq:soln-Z}
\end{equation}
where
\[
Y_1(r)=\int_0^r \frac{KV(\rho)}{\rho}\,d\rho.
\]
Recall that $|\alpha(r)-1|\leq c\,\Om(r)$ 
and $V$ satisfies $M_p(V,r)\leq c\,\Om(r)$ as $r\to 0$.
Moreover,
\[
M_p\left[\int_0^r\frac{KV(\rho)}{\rho}\,d\rho,r\right]\leq c\int_0^r \frac{\Om^2(\rho)}{\rho}\,d\rho,
\]
so all the terms after the $Y(0)$ inside the parentheses of (\ref{eq:soln-Z}) are bounded in $M_p$ either by $\Om(r)$ or by $\int_0^r \Om^2(\rho)\rho^{-1}\,d\rho$ for $0<r<1$. Let us explore $E(r)$.
 Notice that
\[
\int_r^1\frac{Q(\rho)}{\rho}d\rho=\int_r^1[-\alpha_0(\rho)+n\alpha(\rho)]\frac{d\rho}{\rho}
+\int_r^1 Q(\rho)[1-\alpha(\rho)]\frac{d\rho}{\rho}
\]
\[
=\frac{1}{|\del B_1|}\int_{B_1\backslash B_{r}}(na_{ij}(y)y_iy_j|y|^{-2}-a_{ii}(y))\frac{dy}{|y|^n}
+
\int_r^1 Q(\rho)[1-\alpha(\rho)]\frac{d\rho}{\rho}.
\]
But $|Q(\rho)[1-\alpha(\rho)]|\leq c\, \Omega^2(\rho)$, so
\[
\int_r^1Q(\rho)[1-\alpha(\rho)]\,\frac{d\rho}{\rho}=
C-\int_0^rQ(\rho)[1-\alpha(\rho)]|\,\frac{d\rho}{\rho}.
\]
Exponentiating, we obtain
\[
E(r)=\exp\left[\int_r^1\frac{Q(\rho)}{\rho}\,d\rho\right]
=C\exp\left[\frac{1}{|\del B_1|}\int_{B_1\backslash B_{r}}(na_{ij}(y)y_iy_j|y|^{-2}-a_{ii}(y))\frac{dy}{|y|^n}\right]
(1+\zeta_0(r)),
\]
where $|\zeta_0(r)|\leq c\int_0^r\Om^2(\rho)\rho^{-1}\,d\rho$. Since we can rescale $Z$ to set $CY(0)=1$, we
have the formula (\ref{eq:Z-asym}).
This completes the proof of Theorem 1.

\begin{remark}
Choosing $C=-|\del B_1|^{-1}$ in the proof of Theorem 1 and modifying the rest of the argument, we could obtain the asymptotic representation of the fundamental solution of equation (\ref{eq:pde}):
\[
\Gamma_{\mathcal A}(x)=\Gamma(x)+\Lap^{-1}\del_i\del_j((a_{ij}-\de_{ij})\Gamma) +v(x),
\]
where $\Gamma$ is the fundamental solution of $\Lap$ and
\[
M_p(v,r)\leq c\,\Om^2(r)\,r^{2-n}\quad\hbox{for}\ r\in(0,1).
\]
Here the condition (\ref{eq:Om1}) is not necessary; we only need the smallness of $\Om(r)$. However, this formula, as well as the additional terms in the asymptotic expansion of $\Gamma_{\mathcal A}$ can be obtained more easily by iteration from the equation ${\mathcal A}u(x)=\de(x)$.
\end{remark}

\section{Proof of Theorem 2}

\bigskip
Let $q=p/(p-1)$ and $\beta,\gamma\in \RR$. Let us introduce the weighted $L^p$ norm for functions on $\RR^n$ with separate weights at the origin and infinity
\begin{equation}
\|u\|^q_{L^q_{\beta,\gamma}(\RR^n)}=\|u\|^q_{L^q_{\beta}(B_1)}+
\|u\|^q_{L^q_{\gamma}(B_1^c)}
=\int_{|x|<1}|u(x)|^q\,|x|^{\beta q}\,dx +\int_{|x|>1}|u(x)|^q\,|x|^{\gamma q}\,dx,
\label{eq:weightedLq}
\end{equation}
and the weighted Sobolev space $W^{2,q}_{\beta,\gamma}(\RR^n\backslash\{0\})$ with norm
\begin{equation}
\sum_{|\alpha|\leq 2} \|r^{|\alpha|}\partial^\alpha\, u\|_{L^q_{\beta,\gamma}(\RR^n)}.
\label{eq:weightedSobolev}
\end{equation}
Notice that $L^p_{-\beta,-\gamma}(\RR^n)$ is the dual space for $L^q_{\beta,\gamma}(\RR^n)$ and the notation $W^{-2,p}_{-\beta,-\gamma}(\RR^n\backslash\{0\})$ will be used for the dual of $W^{2,q}_{\beta,\gamma}(\RR^n\backslash\{0\})$. 
Many authors have used similar weighted Sobolev spaces to study operators like the Laplacian
on $\RR^n$, $\RR^n\backslash\{0\}$, and other noncompact manifolds with conical or cylindrical ends.

Using the analysis in \cite{MP},  \cite{M} or \cite{LM}, for example,
it is easily verified that the bounded operator
\begin{equation}
\Lap: W^{2,q}_{\beta,\gamma}(\RR^n\backslash\{0\})\to L^q_{\beta+2,\gamma+2}(\RR^n)
\label{eq:Lap-weighted}
\end{equation}
is Fredholm (finite nullity and finite deficiency) for all values of $\beta$ and $\gamma$ {\it except}
for the values $-2+\frac{n}{p}+k$ and $-\frac{n}{q}-k$ where $k$ is any nonnegative integer. In fact, 
(\ref{eq:Lap-weighted}) is an isomorphism for $-n/q<\beta,\gamma<-2+n/p$ (recall that we are assuming $n\geq 3$, so such $\beta,\gamma$ exist). Since we are principally interested in the behavior of functions at the origin, we will fix $\gamma_0\in (-n/q,-2+n/p)$. Then,
for $\beta\in(-2+n/p,-1+n/p)$, we find that (\ref{eq:Lap-weighted}) is surjective with a one-dimensional nullspace
spanned by $|x|^{2-n}$.

Now let us consider the formal adjoint of $\mathcal A$, which also defines a bounded operator on these spaces
\begin{equation}
{\mathcal L}=\overline{a}_{ij}(x)\partial _i\partial _j :  W^{2,q}_{\beta,\gamma} (\RR^n\backslash\{0\}) \to L_{\beta+2,\gamma+2}(\RR^n),
\label{eq:A*-weighted}
\end{equation} 
where $\overline{a}_{ij}$, of course, denotes the complex conjugate of ${a}_{ij}$.
Because $a_{ij}(x)=\de_{ij}$ for $|x|>1$ and $a_{ij}(x)-\de_{ij}$ vanishes as $x\to 0$, the analysis in the above references shows that the operator (\ref{eq:A*-weighted}) is Fredholm for exactly the same values of $\beta$ and $\gamma$ as for 
(\ref{eq:Lap-weighted}).
In fact, for fixed nonexceptional values of $\beta$ and $\gamma$, we may take $\de$ sufficiently small and use perturbation theory (cf.\ \cite{K}, Ch.IV, Sec.5) to conclude that the nullity and deficiency of 
(\ref{eq:Lap-weighted}) and (\ref{eq:A*-weighted}) agree.

So, in addition to the fixed $\gamma_0\in (-n/q,-2+n/p)$, let us now fix $\beta_1 \in (-n/q, -2+n/p)$ and $\beta_2\in(-2+n/p,-1+n/p)$, and denote the adjoints of the corresponding operators (\ref{eq:A*-weighted}) by ${\mathcal A}_1$ and ${\mathcal A}_2$:
\begin{equation}
{\mathcal A}_1 : L^p_{-\beta_1-2,-\gamma_0-2}(\RR^n) \to W_{-\beta_1,-\gamma_0}^{-2,p} (\RR^n\backslash\{0\})
\label{eq:A1}
\end{equation}
is an isomorphism, and 
\begin{equation}
{\mathcal A}_2 : L^p_{-\beta_2-2,-\gamma_0-2}(\RR^n) \to W_{-\beta_2,-\gamma_0}^{-2,p} (\RR^n\backslash\{0\})
\label{eq:A2}
\end{equation}
is injective  with a one-dimensional cokernel. An arbitrary non-zero functional in $\hbox{Coker} \, {\mathcal A}_2$ will be denoted by $\zeta$.  

  We introduce a cut-off function $\eta\in C_0^\infty(B_1)$ equal to $1$ on $B_{1/2}$.
  It follows from (\ref{eq:soln-Z}) that $\eta Z\in L^p_{-\beta_1-2}(B_1)$ but it is not in 
  $L^p_{-\beta_2-2}(B_1)$.
Let $F = {\mathcal A}( (1-\eta) Z)=-{\mathcal A}( \eta Z)$. Since $F = 0$ on $B_{1/2}$ and on $B_1^c$, it follows that $F\in W^{-2,p}_{-\beta_2,-\gamma_0}(\RR^n\backslash\{0\}) \subset 
W^{-2,p}_{-\beta_1,-\gamma_0}(\RR^n\backslash\{0\})$. But (\ref{eq:A1}) is an isomorphism, so $\eta Z$ is the only  solution of ${\mathcal A}_1 v = -F$. Since $L^p_{-\beta_2-2} (B_1)\subset L^p_{-\beta_1-2} (B_1)$, ${\mathcal A}_2 v = -F$ has no solutions,  which means that $\zeta(F)\neq 0$. 

Now let $u\in L^p_{\loc}(B_1)$ be a weak solution of ${\mathcal A}u=0$ satisfying $M_p(u,r)\leq c\,r^{2-n+\e_0}$.
This estimate on $u$ implies that $u\in L^p_{-\beta-2}(B_1)$ for $\beta<-n/q+\e_0$. So let us  restrict our choice of $\beta_1$ to the interval $(-n/q,-n/q+\e_0)$.
Denote by $C$ an arbitrary constant. Since
${\mathcal A} (\eta(u - CZ)) = 0$ on $B_{1/2}$ and on $B_1^c$, we have
${\mathcal A} (\eta(u - CZ))\in W^{-2,p}_{-\beta_2-2,-\gamma_0-2}(\RR^n\backslash\{0\})$.
Choosing $C$ to satisfy
$$C\zeta(F) = \zeta ({\mathcal A} (\eta u)),$$
we obtain $\zeta({\mathcal A}(\eta (u-CZ)))=0$, which implies $\eta( u-CZ)\in \hbox{Dom} ({\mathcal A}_2)$, and in particular
$$\eta(u- CZ) \in L^p_{-\beta_2-2} (B_1).$$
But this implies that $w=u-CZ$ satisfies  $M_p(w,r)\leq c\,r^{\beta_2+2-n/p}$.
Assuming that we had fixed $\beta_2=\frac{n}{p}-1-\e_1$ where $\e_1\in (0,1)$, we conclude that  
\begin{equation}
M_p(w,r)\leq c\,r^{1-\e_1}.
\label{eq:w-est}
\end{equation}
This completes the proof.


\end{document}